\documentclass[11pt]{amsart}
\usepackage{amssymb,amsmath,txfonts,mathrsfs,mathtools,titletoc, tikz, stmaryrd}
\usepackage[alphabetic]{amsrefs}

\newtheorem{theorem}{Theorem}[section]
\newtheorem{prop}[theorem]{Proposition}
\newtheorem{lemma}[theorem]{Lemma}

\newtheorem{question}[theorem]{Question}

\newtheorem{cor}[theorem]{Corollary}

\numberwithin{equation}{section}

\def\pf{{\it Proof:}~}

\begin{document}

\title[Integral of scalar curvature]{Integral of scalar curvature on non-parabolic manifolds}
\author{Guoyi Xu}
\address{Department of Mathematical Sciences\\Tsinghua University, Beijing\\P. R. China, 100084}
\email{guoyixu@tsinghua.edu.cn}
\date{\today}
\date{\today}

\begin{abstract}
Using the monotonicity formulas of Colding and Minicozzi, we prove that on any complete, non-parabolic Riemannian manifold $(M^3, g)$ with non-negative Ricci curvature, the asymptotic weighted scaling invariant integral of scalar curvature has an explicit bound in form of  asymptotic volume ratio.
\\[3mm]
\end{abstract}
\thanks{The author was partially supported by NSFC 11771230}

\maketitle

\section{Introduction}

The study of integral of the curvature started from the well-known Gauss-Bonnet Theorem: for any compact $2$-dim Riemannian manifold $(M^2, g)$, $\int_{M^2} K d\mu= 2\pi\chi(M^2)$, where $K$ is the Gaussian curvature of $M^2$, $d\mu$ is the element of area of $M^2$, $\chi(M^2)$ is the Euler characteristic of $M^2$.

Cohn-Vossen \cite{CV} studied the integral of the curvature on complete $2$-dim Riemannian manifold, obtained the so-called Cohn-Vossen's inequality: If $(M^2, g)$ is a finitely connected, complete, oriented Riemannian manifold, and assume $\int_{M^2} Kd\mu$ exists as extended real number, then $\int_{M^2} Kd\mu\leq 2\pi\chi(M^2)$.

Motivated to get a generalization of the Cohn-Vossen's inequality, Yau \cite[Problem $9$]{Yau} posed the following question: Given a $n$-dimensional complete manifold $(M^n, g)$ with $Rc\geq 0$, let $B_r(p)$ be the geodesic ball around $p\in M^n$ and $\sigma_k$ be the $k$-th elementary symmetric function of the Ricci tensor, is it true that $\lim\limits_{r\rightarrow \infty} r^{-n+2k} \int_{B_r(p)} \sigma_k< \infty$?

In $2013$, Bo Yang \cite{YangBo} constructed examples, which answered the above question for $k> 1$ negatively. However, the interesting case $k= 1$ is still open, where $\sigma_1$ is the scalar curvature $R$. We formulate it in the following question separately. 
\begin{question}[Yau]\label{ques Yau-1}
{For any complete Riemannian manifold $(M^n, g)$ with $Rc\geq 0$, any $p\in M^n$, is it true that $\lim\limits_{r\rightarrow \infty} r^{2- n}\int_{B_r(p)} R< \infty$ ?
}
\end{question}

Related to the above question, Shi and Yau \cite{ShiYau} gave a scaling invariant upper bound estimate for the average integral of the scalar curvature, on K$\ddot{a}$hler manifolds with bounded, pinched, nonnegative holomorphic bisectional curvature.

On the other hand, if we relax the assumption $Rc\geq 0$ to the non-negative sectional curvature $K\geq 0$, among other things Petrunin \cite{Petrunin} proved: There exists $C(n)> 0$, such that for any complete Riemannian manifold $(M^n, g)$ with sectional curvature $K\geq 0$ and any $p\in M^n$, $\int_{B_1(p)} R\leq C(n)$. Petrunin's result implies that Question \ref{ques Yau-1} has one partial affirmative answer when $K\geq 0$.

A complete Riemannian manifold is said to be \textbf{non-parabolic} if it admits a positive Green function, otherwise it is said to be \textbf{parabolic}. By a result of Varopoulos \cite{Var} a complete non-compact Riemannian manifold with $Rc\geq 0$ is non-parabolic if and only if 
\begin{align}
\int_1^\infty \frac{r}{\mathrm{Vol}(B_r(p))}dr< \infty \nonumber 
\end{align}

For non-parabolic Riemannian manifolds, there exists a unique, minimal, positive Green function, denoted as $G(p, x)$, where $p$ is a fixed point on manifold. We define $b(x)= [n(n- 2)\omega_n\cdot G(p, x)]^{\frac{1}{2- n}}$, where $\omega_n$ is the volume of the unit ball in $\mathbb{R}^n$. We also define the \textbf{asymptotic volume ratio} of the manifold $M^n$ as: $\mathrm{V}_{M}= \lim\limits_{r\rightarrow \infty}\frac{\mathrm{Vol}(B_r(p))}{\omega_nr^n}$. 

In this short note, we proved the following theorem, which makes some progress to Question \ref{ques Yau-1} in $3$-dim case. 
\begin{theorem}\label{thm 3-dim upper bound of av integ cur}
{For a complete non-compact Riemannian manifold $M^3$, which is non-parabolic with $Rc\geq 0$, we have 
\begin{align}
\varlimsup_{r\rightarrow \infty}\frac{\int_{b\leq r}R\cdot |\nabla b|}{r}\leq 8\pi \big[1- \mathrm{V}_{M}\big] \nonumber 
\end{align}
}
\end{theorem}

\section{The estimates of Green function and its relatives}

In the rest of the paper, we always use $\rho(x)= d(p, x)$ unless otherwise mentioned.  From the behavior of Green function $G(p, x)$ near singular point $p$, we get 
\begin{align}
\lim_{\rho(x)\rightarrow 0}\frac{b(x)}{\rho(x)}&= 1 \label{estimate of b(0)} \\
\lim_{\rho(x)\rightarrow 0} |\nabla b|&= 1 \label{est of gradient of b at 0}
\end{align}

Define $\theta_p(r)= \frac{\mathrm{Vol}(\partial B_r(p))}{r^{n- 1}}$ and $\theta= \lim\limits_{r\rightarrow \infty}\theta_p(r)$, from Bishop-Gromov Volume Comparison Theorem, the limit always exists. 

The following Lemma was essentially proved in \cite{CM-AJM} firstly, which used Gromov-Hausdorff convergence. Our statement followed from the intrinsic argument in \cite{LTW}.

\begin{lemma}\label{lem point-wise bound of b}
{If $M^n$ has $Rc\geq 0$ with $n\geq 3$ and maximal volume growth, for any $\delta \in (0, \frac{1}{2}]$, we have 
\begin{align}
\big(\mathrm{V}_M\big)^{\frac{1}{n- 2}}\Big(1+ \tau\Big)^{\frac{1}{2- n}}\rho(x) \leq b(x)\leq \big(\mathrm{V}_M\big)^{\frac{1}{n- 2}}\Big(1- \tau\Big)^{\frac{1}{2- n}}\rho(x) \nonumber 
\end{align}
where $\tau= C(n)\big[\delta+ (\theta_p(\delta \rho(x))- \theta)^{\frac{1}{n- 1}}\big]$. Especially, $\lim\limits_{\rho(x)\rightarrow \infty} \frac{b(x)}{\rho(x)}= \big(\mathrm{V}_M\big)^{\frac{1}{n- 2}}$.
}
\end{lemma}\qed

\begin{lemma}\label{lem upper bound of gradient of b}
{If $M^n$ has $Rc\geq 0$ with $n\geq 3$, and it is non-parabolic and not maximal volume growth, then $\lim\limits_{r\rightarrow \infty} \sup\limits_{\rho(x)= r}|\nabla b|(x)= 0$. Furthermore, if it has maximal volume growth, then $|\nabla b|\leq 1$.
}
\end{lemma}

\pf
{By Li-Yau's lower bound for the Green function \cite{LY}, 
\begin{align}
C\int_{\rho(x)}^{\infty} \frac{s}{\mathrm{Vol}(B_s(p))}ds\leq G(p, x) \nonumber 
\end{align}

Then for $\rho(x)= r_0$, apply Bishop-Gromov Volume Comparison Theorem, 
\begin{align}
G(p, x)\geq \frac{C}{r_0^{-n}\mathrm{Vol}(B_{r_0}(p))}\rho^{2- n}(x) \nonumber 
\end{align}

By Cheng-Yau \cite{CY} gradient estimate at such $x$, 
\begin{align}
|\nabla b|(x)= C(n)G^{\frac{1}{2- n}}(p, x)\cdot \big|\frac{\nabla G}{G}\big|(p, x)\leq C(n)\frac{G^{\frac{1}{2- n}}(p, x)}{\rho(x)}\leq C(n)\Big(r_0^{-n}\mathrm{Vol}\big(B_{r_0}(p)\big)\Big)^{\frac{1}{n- 2}} \nonumber 
\end{align} 
If $M^n$ has not maximal volume growth, we have the conclusion.

If $M^n$ has maximal volume growth, the conclusion follows from \cite[Theorem $3.1$]{Colding}.
}
\qed

The following lemma was implied by the argument in \cite{CC-Ann}, and was used repeatedly in \cite{CM-AJM}, \cite{Colding} and \cite{CM}. We give a direct proof of this result here for reader's convienence.
\begin{lemma}\label{lem integral gradient est of b}
{If $M^n$ has $Rc\geq 0$ with $n\geq 3$ and maximal volume growth,  we have $\lim\limits_{r\rightarrow \infty}\frac{\int_{\hat{b}\leq r} |\nabla \hat{b}- \nabla \rho|^2}{V(\hat{b}\leq r)}= 0$, where $\hat{b}= \big(\mathrm{V}_M\big)^{\frac{1}{2- n}}\cdot b$.
}
\end{lemma}

\pf
{Firstly we recall that $\int_{b\leq r}|\nabla b|^2= \omega_nr^n$, which implies
\begin{align}
\int_{\hat{b}\leq r}|\nabla \hat{b}|^2= \mathrm{V}_M\omega_n r^n \nonumber 
\end{align}

From the Green's formula and $\Delta b= (n- 1)\frac{|\nabla b|^2}{b}$, we get
\begin{align}
\int_{\hat{b}\leq r}\nabla \hat{b}\cdot \nabla \rho&= -\int_{\hat{b}\leq r}\Delta \hat{b}\cdot \rho+ \int_{\hat{b}= r}|\nabla \hat{b}|\cdot \rho \nonumber \\
& = -(n- 1)\int_{\hat{b}\leq r}\rho\frac{|\nabla \hat{b}|^2}{\hat{b}}+ \int_{\hat{b}= r}|\nabla \hat{b}|r+ \int_{\hat{b}= r}|\nabla \hat{b}|\cdot (\rho- r) \nonumber 
\end{align}
use $\int_{b= r}|\nabla b|= n\omega_n r^{n- 1}$, we have
\begin{align}
\int_{\hat{b}= r}|\nabla \hat{b}|r= \mathrm{V}_Mn\omega_n r^{n} \nonumber 
\end{align}

Note by Lemma \ref{lem point-wise bound of b} and Lemma \ref{lem upper bound of gradient of b}, the following holds
\begin{align}
\lim_{r\rightarrow\infty} \frac{\int_{\hat{b}\leq r}\rho\frac{|\nabla \hat{b}|^2}{\hat{b}}- |\nabla \hat{b}|^2}{V(\hat{b}\leq r)}= 0 \quad \quad and \quad \quad \lim_{r\rightarrow\infty} \frac{\int_{\hat{b}= r}|\nabla \hat{b}|\cdot (\rho- r)}{V(\hat{b}\leq r)}= 0 \nonumber 
\end{align}

By the above, using Lemma \ref{lem point-wise bound of b} again,
\begin{align}
\lim_{r\rightarrow\infty} \frac{\int_{\hat{b}\leq r}\nabla \hat{b}\cdot \nabla \rho}{V(\hat{b}\leq r)}&= -(n- 1)\lim_{r\rightarrow\infty}\frac{\int_{\hat{b}\leq r}|\nabla \hat{b}|^2}{V(\hat{b}\leq r)}+ \lim_{r\rightarrow\infty}\frac{\mathrm{V}_Mn\omega_nr^n}{V(\hat{b}\leq r)} \nonumber \\
&= \lim_{r\rightarrow\infty}\frac{\mathrm{V}_M\omega_nr^n}{V(\hat{b}\leq r)}= 1\nonumber
\end{align}

Finally, 
\begin{align}
\lim\limits_{r\rightarrow \infty}\frac{\int_{\hat{b}\leq r} |\nabla \hat{b}- \nabla \rho|^2}{V(\hat{b}\leq r)}&= \lim\limits_{r\rightarrow \infty}\frac{\int_{\hat{b}\leq r} |\nabla \hat{b}|^2}{V(\hat{b}\leq r)}+ \lim\limits_{r\rightarrow \infty}\frac{\int_{\hat{b}\leq r} |\nabla \rho|^2}{V(\hat{b}\leq r)}- 2\lim_{r\rightarrow\infty} \frac{\int_{\hat{b}\leq r}\nabla \hat{b}\cdot \nabla \rho}{V(\hat{b}\leq r)} \nonumber \\
&= \lim\limits_{r\rightarrow \infty}\frac{\mathrm{V}_M\omega_nr^n}{V(\hat{b}\leq r)}+ 1- 2= 0 \nonumber 
\end{align}

}
\qed

\begin{cor}\label{cor limit of the rescaled volume by b}
{For a complete non-compact Riemannian manifold $M^n$, which is non-parabolic with $Rc\geq 0$, we have $\lim\limits_{r\rightarrow \infty}\frac{\int_{b\leq r}|\nabla b|^3}{r^n}= \big(\mathrm{V}_M\big)^{\frac{1}{n- 2}}\omega_n$.
}
\end{cor}

\pf
{If $M^n$ does not have the maximal volume growth, from Lemma \ref{lem upper bound of gradient of b}, the conclusion follows directly.

In the rest of the proof, we assume that $M^n$ has maximal volume growth. Let $\tilde{r}= \big(\mathrm{V}_M\big)^{\frac{1}{2- n}}r, \hat{b}= \big(\mathrm{V}_M\big)^{\frac{1}{2- n}}b$, we have 
\begin{align}
\frac{\int_{b\leq r}\Big||\nabla b|^3- \big(\mathrm{V}_M\big)^{\frac{3}{n- 2}}\Big|}{r^n}&= \big(\mathrm{V}_M\big)^{\frac{3- n}{n- 2}}\frac{\int_{\hat{b}\leq \tilde{r}}\big||\nabla \hat{b}|^3- 1\big|}{\tilde{r}^n} \nonumber 
\end{align}

From Lemma \ref{lem integral gradient est of b}, Lemma \ref{lem upper bound of gradient of b} and Lemma \ref{lem point-wise bound of b}, use the Bishop-Gromov Volume Comparison Theorem,
\begin{align}
\varlimsup_{r\rightarrow \infty}\frac{\int_{\hat{b}\leq r}\big||\nabla \hat{b}|^3- 1\big|}{r^n}&\leq C(\mathrm{V}_M)\varlimsup_{r\rightarrow \infty}\frac{\int_{\hat{b}\leq r} \big||\nabla \hat{b}|- 1\big|}{r^n}\leq C\cdot \varlimsup_{r\rightarrow \infty}\frac{\int_{\hat{b}\leq r} \big|\nabla \hat{b}- \nabla \rho\big|}{r^n} \nonumber \\
&\leq C\cdot \varlimsup_{r\rightarrow \infty}\frac{\big(\int_{\hat{b}\leq r} \big|\nabla \hat{b}- \nabla \rho\big|^2\big)^{\frac{1}{2}}}{r^n}\cdot V(\hat{b}\leq r)^{\frac{1}{2}} \nonumber \\
&\leq C\cdot \Big(\varlimsup_{r\rightarrow \infty}\frac{\int_{\hat{b}\leq r} \big|\nabla \hat{b}- \nabla \rho\big|^2}{V(\hat{b}\leq r)}\Big)^{\frac{1}{2}}= 0 \nonumber 
\end{align}
which implies $\lim\limits_{r\rightarrow \infty}\frac{\int_{\hat{b}\leq r}\big||\nabla \hat{b}|^3- 1\big|}{r^n}= 0$. 

Hence from the above and Lemma \ref{lem point-wise bound of b}, we have 
\begin{align}
\lim\limits_{r\rightarrow \infty}\frac{\int_{b\leq r}|\nabla b|^3}{r^n}= \big(\mathrm{V}_M\big)^{\frac{3}{n- 2}}\lim\limits_{r\rightarrow \infty}\frac{V(b\leq r)}{r^n}=  \big(\mathrm{V}_M\big)^{\frac{3- n}{n- 2}}\lim\limits_{r\rightarrow \infty}\frac{V(\rho\leq r)}{r^n}= \big(\mathrm{V}_M\big)^{\frac{1}{n- 2}}\omega_n \nonumber 
\end{align}
}
\qed

\section{Integral of the scalar curvature by dimension reduction}

Define 
\begin{align}
\mathcal{A}(r)= r^{1- n}\int_{b= r} |\nabla b|^2 \quad \quad and \quad \quad \mathcal{V}(r)= r^{2- n}\int_{b\leq r}\frac{|\nabla b|^3}{b^2} \nonumber 
\end{align}
from the above definition, it is straightforward to get 
\begin{align}
\mathcal{V}'(r)= r^{-1}(2- n)\mathcal{V}(r)+ r^{-1}\mathcal{A}(r) \label{deri of V-M}
\end{align}

Let $\vec{n}= \frac{\nabla b}{|\nabla b|}$ and $\mathbf{B}= \nabla^2(b^2)- 2|\nabla b|^2\cdot g$, where $g$ is the Riemannian metric on $M^n$. Let $\langle\mathbf{B}(\vec{n}), v\rangle= \mathbf{B}(\vec{n}, v)$ for any $v$, then $|\mathbf{B}(\vec{n})|^2= |\mathbf{B}(\vec{n})^T|^2+ \mathbf{B}(\vec{n},\vec{n})^2$.

The following lemma is the analogue of Theorem $2.12$ of \cite{Colding}.
\begin{lemma}\label{lem limit of A-M and V-M}
{For a complete non-compact Riemannian manifold $M^n$, which is non-parabolic with $Rc\geq 0$, we have 
\begin{align}
\lim_{r\rightarrow \infty}\mathcal{A}(r)= (n- 2)\lim_{r\rightarrow \infty}\mathcal{V}(r)= n\omega_n\big(\mathrm{V}_M\big)^{\frac{1}{n- 2}} \nonumber 
\end{align}
}
\end{lemma}
\pf
{Let $\beta= 1$ in \cite[Theorem $3.2$]{CM}, we obtain
\begin{align}
&\quad r^{1- n}\int_{b\leq r}|\nabla b| \big[|\Pi_0|^2+ Rc(\vec{n}, \vec{n})\big]+ \frac{|\mathbf{B}(\vec{n})|^2+ (n- 2)|\mathbf{B}(\vec{n})^T|^2}{4(n- 1)b^2|\nabla b|} \nonumber \\
&=  \mathcal{A}'(r)- 2(n- 2)\mathcal{V}'(r) \label{upper bound of integral along direction of b}
\end{align}

From \cite[Theorem $1.1$]{CM}, we have  
\begin{align}
\mathcal{A}'(r)\leq 0 \label{deri of A is negative}
\end{align}
Combining (\ref{upper bound of integral along direction of b}), we get
\begin{align}
\mathcal{V}'(r)\leq \mathcal{A}'(r)\leq 0 \label{deri of A-M is non-positive}
\end{align}

From $\mathcal{A}(r), \mathcal{V}(r)\geq 0$ and (\ref{deri of A-M is non-positive}), we know that $\lim\limits_{r\rightarrow \infty}\mathcal{A}(r), \lim\limits_{r\rightarrow \infty}\mathcal{V}(r)$ exists. By L'H\^{o}pital's rule,
\begin{align}
\lim_{r\rightarrow \infty}\mathcal{V}(r)= \lim_{r\rightarrow \infty}\frac{\int_{b\leq r} b^{-2}|\nabla b|^3}{r^{n- 2}}= \lim_{r\rightarrow \infty}\frac{\int_{b= r}b^{-2}|\nabla b|^2}{(n- 2)r^{n- 3}}= \lim_{r\rightarrow \infty}\frac{1}{n- 2}\mathcal{A}(r) \nonumber 
\end{align}

On the other hand, from the L'H\^{o}pital's rule,
\begin{align}
\lim\limits_{r\rightarrow \infty}\frac{\int_{b\leq r}|\nabla b|^3}{r^n}= \lim_{r\rightarrow \infty}\frac{\int_{b= r}|\nabla b|^2}{nr^{n- 1}}= \lim_{r\rightarrow \infty}\frac{1}{n}\mathcal{A}(r) \label{limit of A-para}
\end{align}

The conclusion follows from Corollary \ref{cor limit of the rescaled volume by b} and (\ref{limit of A-para}).
}
\qed



\begin{prop}\label{prop integral along b direction}
{For a complete non-compact Riemannian manifold $M^n$, which is non-parabolic with $Rc\geq 0$, we have
\begin{align}
&\lim_{r\rightarrow \infty}\frac{\int_{b\leq r}|\nabla b|\big[|\Pi_0|^2+ Rc(\vec{n})\big]+ \frac{|\mathbf{B}(\vec{n})|^2+ (n- 2)|\mathbf{B}(\vec{n})^T|^2}{4(n- 1)b^2|\nabla b|}}{r^{n- 2}}= 0 \nonumber \\
&\lim_{r\rightarrow \infty}r^{2- n}\int_{b\leq r}|\nabla b|\cdot H^2= \frac{(n- 1)^2n}{n- 2}\omega_n\big(\mathrm{V}_M\big)^{\frac{1}{n- 2}}  \nonumber 
\end{align}
where $H$ is the mean curvature of the level set of $b$ with respect to the normal vector $\frac{\nabla b}{|\nabla b|}$.
}
\end{prop}

\pf
{From (\ref{upper bound of integral along direction of b}), (\ref{deri of A is negative}) and (\ref{deri of V-M}), we have 
\begin{align}
&\quad \varlimsup_{r\rightarrow \infty}\frac{\int_{b\leq r}|\nabla b|\big[|\Pi_0|^2+ Rc(\vec{n})\big]+ \frac{|\mathbf{B}(\vec{n})|^2+ (n- 2)|\mathbf{B}(\vec{n})^T|^2}{4(n- 1)b^2|\nabla b|}}{r^{n- 2}} \nonumber \\
&= \varlimsup_{r\rightarrow \infty}r\Big[\mathcal{A}'(r)- 2(n- 2)\mathcal{V}'(r)\Big] \nonumber \\
&\leq \varlimsup_{r\rightarrow \infty} 2(n- 2)\big[(n- 2)\mathcal{V}(r)- \mathcal{A}(r)\big]= 0 \label{ineq 1.1 need} 
\end{align}
the last equation above follows from Lemma \ref{lem limit of A-M and V-M}.

It is straightforward to compute the mean curvature of the level set of $b$ with respect to the normal vector $\frac{\nabla b}{|\nabla b|}$ as the following:
\begin{align}
H= \frac{(n- 1)|\nabla b|}{b}- \frac{\mathbf{B}(\vec{n}, \vec{n})}{2b|\nabla b|} \label{expression of H in form of b}
\end{align}

From (\ref{expression of H in form of b}) and (\ref{ineq 1.1 need}), also note $\mathbf{B}(\vec{n}, \vec{n})^2\leq |\mathbf{B}(\vec{n})|^2$, we get
\begin{align}
\lim_{r\rightarrow \infty}r^{2- n}\int_{b\leq r}|\nabla b|\cdot H^2&= \lim_{r\rightarrow \infty}r^{2- n}\int_{b\leq r} \Big\{(n- 1)^2\frac{|\nabla b|^3}{b^2}- (n- 1)\frac{|\nabla b|\cdot \mathbf{B}(\vec{n}, \vec{n})}{b^2}+ \frac{|\mathbf{B}(\vec{n}, \vec{n})|^2}{4b^2|\nabla b|}\Big\} \nonumber \\
&\leq (1+ \epsilon)\lim_{r\rightarrow \infty}\frac{\int_{b\leq r} (n- 1)^2\frac{|\nabla b|^3}{b^2}}{r^{2- n}}+ \big(\frac{4}{\epsilon}+ 1\big)\lim_{r\rightarrow \infty}
\frac{\int_{b\leq r} \frac{|\mathbf{B}(\vec{n})|^2}{4b^2|\nabla b|}}{r^{n- 2}} \nonumber \\
&\leq (1+ \epsilon)\lim_{r\rightarrow \infty}\frac{\int_{b\leq r} (n- 1)^2\frac{|\nabla b|^3}{b^2}}{r^{2- n}} \nonumber 
\end{align}
let $\epsilon\rightarrow 0$ in the above, we have 
\begin{align}
\lim_{r\rightarrow \infty}r^{2- n}\int_{b\leq r}|\nabla b|\cdot H^2\leq \lim_{r\rightarrow \infty}\frac{\int_{b\leq r} (n- 1)^2\frac{|\nabla b|^3}{b^2}}{r^{2- n}} = (n- 1)^2\lim_{r\rightarrow \infty}\mathcal{V}(r) \nonumber 
\end{align}

Similar as the above, we can get $\lim\limits_{r\rightarrow \infty}r^{2- n}\int_{b\leq r}|\nabla b|\cdot H^2\geq (n- 1)^2\lim\limits_{r\rightarrow \infty}\mathcal{V}(r)$. The conclusion follows from the above argument and Lemma \ref{lem limit of A-M and V-M}.
}
\qed

\begin{prop}\label{prop dim reduction method}
{For a complete non-compact Riemannian manifold $M^n$, which is non-parabolic with $Rc\geq 0$, we have
\begin{align}
\lim_{r\rightarrow \infty}\frac{\int_{b\leq r}|\nabla b|R- \int_0^rdt\int_{b= t} R\big(b^{-1}(t)\big)}{r^{n- 2}}= -(n- 1)n\omega_n\big(\mathrm{V}_M\big)^{\frac{1}{n- 2}}  \nonumber 
\end{align}
}
\end{prop}

\pf
{From the Gauss equation, we have 
\begin{align}
R(M^n)&= R\big(b^{-1}(t)\big)+ 2Rc(\vec{n})- 2\sum_{i\neq j}\lambda_i\lambda_j \nonumber \\
&= R\big(b^{-1}(t)\big)+ 2Rc(\vec{n})- \Big(\frac{n- 2}{n- 1}H^2- |\Pi_0|^2\Big) \nonumber \\
&= R\big(b^{-1}(t)\big)+ 2Rc(\vec{n})+ |\Pi_0|^2- \frac{n- 2}{n- 1}H^2 \nonumber 
\end{align}

The conclusion follows from Proposition \ref{prop integral along b direction}.
}
\qed

In the rest of this section, we will apply the general results obtained before, to study the curvature behavior on $3$-dim Riemannian manifolds and their applications.

\begin{lemma}\label{lem 3-dim level set is connected}
{For a complete non-compact Riemannian manifold $M^3$, which is non-parabolic and diffeomorphic to $\mathbb{R}^3$, if $b^{-1}(t)$ is a smooth surface, then it is connected.
}
\end{lemma}

\pf
{If $\Omega_1, \Omega_2$ are two connected components of $b^{-1}(t)$, then the base point of $b$, denoted as $p$, is enclosed by the unique surface $\Omega_i$ (otherwise $G\equiv t^{-1}$ on the region enclosed by $\Omega_1$ and $\Omega_2$ from the maximum principle, from the unique continuation of harmonic function, we get the contradiction). 

Note $M^3$ is diffeomorphic to $\mathbb{R}^3$, hence the other $\Omega_{i+ 1}$ encloses one region $\Omega\subset M^3$, and $p\notin \Omega$. By the maximum principle again $G\equiv t^{-1}$ on $\Omega$, the contradiction follows from the unique continuation of harmonic function again.

}
\qed

\pf [The proof of Theorem \ref{thm 3-dim upper bound of av integ cur}]
{If the universal cover $\tilde{M}$ of $M^3$ is isometric to $N^2\times \mathbb{R}$, where $K(N^2)\geq 0$; because $M^3$ is non-parabolic, we know that $N^2$ is complete and non-compact. If $K(N)\equiv 0$, then 
\begin{align}
\lim_{s\rightarrow \infty} s^{-1}\int_{B_s(p)} R(M) d\mu= 0 \nonumber 
\end{align}

Otherwise, $K(N)> 0$ at some point, from \cite{Soul-P}, we know that $N^2$ is diffeomorphic to $\mathbb{R}^2$, hence the universal cover $\tilde{M}$ is diffeomorphic to $\mathbb{R}^3$. If the universal cover of $M^3$ is not isometric to $N^2\times \mathbb{R}$, then from \cite[Theorem $2$]{LiuGang}, $M^3$ is also diffeomorphic to $\mathbb{R}^3$. 

In the above two cases, from the Gauss-Bonnet Theorem and Lemma \ref{lem 3-dim level set is connected}, we have
\begin{align}
\int_0^rdt\int_{b= t} R\big(b^{-1}(t)\big)\leq 8\pi\cdot r \label{integral of level set curv}
\end{align}
then conclusion follows from (\ref{integral of level set curv}) and Proposition \ref{prop dim reduction method}.
}
\qed

\section*{Acknowledgments}
The author thank the anonymous referee for helpful suggestion, especially providing a simple proof of Lemma \ref{lem integral gradient est of b}.

\end{document}